\documentclass[12pt,a4paper]{article}
\usepackage{amsthm}
\usepackage[cp1251]{inputenc}   % Windows
\usepackage[english, russian]{babel}

\usepackage{amsmath}
\usepackage{amsfonts,amssymb}

\usepackage{graphicx}                 % Для включения графических изображений

\usepackage{psfrag}                    % Надписи на рисунке средствами \LaTeX'а
\usepackage{colordvi}

\usepackage{cite}               % Для корректного оформления литературы
\usepackage{srcltx}            % Для удобного перехода между WiNEdT и YAP
\usepackage{url}
\usepackage{amsxtra}
\usepackage{algorithm}
\usepackage{algorithmic}
\usepackage{caption}
\usepackage{makeidx}
\usepackage{multirow}

\usepackage{mathtools}
\usepackage{algorithm}
\DeclarePairedDelimiter\bracket{(}{)}
\newcommand{\br}[1]{\bracket*{#1}}
\DeclarePairedDelimiter\figbracket{\{}{\}}
\newcommand{\fbr}[1]{\figbracket*{#1}}

\newtheorem{example}{Пример}
\newtheorem{remark}{Замечание}

\newtheorem{theorem}{Теорема}

\usepackage{enumerate}

\newcommand\abs[1]{\left|#1\right|}

\def\ge{\geqslant}
\def\le{\leqslant}
\def\R{{\mathbb R}}

\begin{document}
\begin{center}
{\bf ОБ ОДНОМ МЕТОДЕ МИНИМИЗАЦИИ ВЫПУКЛОЙ ЛИПШИЦЕВОЙ ФУНКЦИИ ДВУХ ПЕРЕМЕННЫХ НА КВАДРАТЕ\footnote{Исследования Ф. С. Стонякина, посвящённые схеме доказательства теоремы 1, а также разработке раздела 5.2, а также замечаний 2, 5, 6 и примера 2 выполнены при поддержке гранта РНФ 18-71-00048.}
\\[2ex]
Д.~А.~Пасечнюк, Ф.~С.~Стонякин}\\[3ex]
\end{center}

%%%%%%%%%%%%%%%%%%%%%%%%%%%%%%
%%%%%%%%%%%%%%%%%%%%%%%%%%%%%%
%%%%%%ВВЕДЕНИЕ%%%%%%%%%%%%%%%%
%%%%%%%%%%%%%%%%%%%%%%%%%%%%%%
%%%%%%%%%%%%%%%%%%%%%%%%%%%%%%
\section{Введение}

В последние десятилетия активно исследуются самые разные подходы к численным методам для условных и безусловных задач многомерной оптимизации \cite{vasiliev2011methods, gasnikov2017book, bubeck2015convex, nesterov2018lect}. В частности, условные задачи оптимизации с небольшим количеством функциональных ограничений могут с помощью метода множителей Лагранжа могут сводиться к маломерным двойственным задачам. А для задач маломерной оптимизации обычно используют разные геометрические методы, среди которых весьма популярны метод центров тяжести Левина-Ньюмена, а также метод эллипсоидов (см., например, \cite{gasnikov2017book, bubeck2015convex, nesterov2018lect}). Отметим \cite{nemirovcky1979complexity, bubeck2015convex}, что метод центров тяжести имеет оптимальную оценку числа итераций $n$ на классе выпуклых функционалов для достижения необходимой точности решения $\varepsilon$ задачи ($d$ --- размерность пространства):
\begin{equation}\label{LinEstimate}
n \geqslant Cd\log\left(\frac{1}{\varepsilon}\right).
\end{equation}
Однако стоимость самих итераций метода центра тяжести может быть довольно дорогой за счёт трудоёмкой операций отыскания центра тяжести выпуклого множества (см., например \cite{gasnikov2017book, bubeck2015convex}). Поэтому на практике часто применяют метод эллипсоидов, где операция нахождения центра тяжести заменяется на задачу нахождения центра описанного эллипсоида, что делает итерацию метода менее затратной. Указанные методы работают с линейной скоростью сходимости (сопоставимую с \eqref{LinEstimate}) как для гладких, так и негладких функций и требуют для достижения достаточно приемлемого качества решения вычисления (суб)градиента на каждой итерации.

Известны самые разные подходы к проблеме ускорения работы методов оптимизации \cite{gasnikov2017book, bubeck2015convex, nesterov2018lect}. Например, для удешевления стоимости итераций метода часто применяют рандомизированные процедуры, предполагающие вместо вычисления точных значений важных параметров (например, градиента) на итерациях их вероятностных аналогов \cite{bayandina2018bezgrad,vorontsova2019accelerated}. Например, вычисление вместо (суб)градиентов на итерациях метода их стохастических аналогов позволяет уменьшить стоимость итерации. Однако для такого типа неполноградиентных методов при этом теряется детерминированные оценки гарантированного качества найденного решения заменяются на вероятностные (стохастические).

Мы же рассмотрим в настоящей статье предложенный Ю.Е. Нестеровым (см., упражнение 1.5 из \cite{gasnikov2017book}) неполноградиентный метод минимизации липшицевой выпуклой липшицевой функции двух переменных на квадрате, позволяющий получить детерминированные оценки.
Особенность метода --- сведение задачи минимизации к постепенному сужению области определения за счёт учёта направления (суб)градиентов целевой функции вблизи точек-решений вспомогательных одномерных задач. Иными словами, предлагается комбинировать быстро работающие методы одномерной минимизации и с учётом геометрической структуры области определения выяснять направление субградиентов в точках на каждой итерации. Важно при этом, что на каждой итерации метода два раза находится направление (суб)градиента (а не самого градиента!), что существенно удешевляет стоимость итерации по сравнению с любым полноградиентным методом. При этом сохраняется линейная скорость сходимости метода.

Основной результат работы --- оценка скорости сходимости метода Ю.Е. Нестерова для гладких выпуклых липшицевых функций, имеющих липшицев градиент (теорема 1). Метод может быть использован и в случае вычисления (направления) градиента с некоторой погрешностью (замечание 2). Отмечено, что свойство липшицевости градиента достаточно потребовать не на всём квадрате, а лишь на некоторых отрезках внутри этого квадрата (замечание 5). Показано, что метод может эффективно работать и для некоторых негладких функций. Однако построен пример простой негладкой функции, для которой при неверном выборе субградиента свойство сходимости метода может не выполняться (пример 1). Проведено сравнение работы метода Ю.Е. Нестерова, метода эллипсоидов и градиентного спуска для гладких выпуклых функций различного вида. Проведённые эксперименты показали, что рассматриваемый метод может позволить достичь требуемой точности решения задачи за меньшее в сравнении с другими рассмотренными методами время как в случае известного точного решения, так и при использовании оценок на скорость сходимости методов. Замечено, что при увеличении точности искомого решения время работы метода Ю.Е. Нестерова может расти медленнее, чем время работы метода эллипсоидов. Описана также методика применения указанного подхода к многомерным задачам выпуклого программирования с двумя функциональными ограничениями.

Отметим, что норму и скалярное произведение на плоскости всюду мы считаем стандартными.

%%%%%%%%%%%%%%%%%%%%%%%%%%%%%%
%%%%%%%%%%%%%%%%%%%%%%%%%%%%%%
%%%%%%ОПИСАНИЕ МЕТОДА%%%%%%%%%
%%%%%%%%%%%%%%%%%%%%%%%%%%%%%%
%%%%%%%%%%%%%%%%%%%%%%%%%%%%%%
\section{Описание метода}
    Дан квадрат $\Pi \subset \mathbb{R}^2$ со стороной $R$, задана функция $f: \Pi \rightarrow \mathbb{R}$.

    Опишем итерацию предлагаемого метода (см. \cite{gasnikov2017book}, упражнение 1.5). Сформулируем его как для гладких, так и негладких целевых функционалов, хотя метод может расходиться на негладкой функции (см. пример 1 далее). При этом отметим, что для негладких выпуклых функционалов вместо обычного градиента используется {\it субградиент} $\upsilon \in \R^2$ определяемый неравенством

\[
    f(x) - f(x_0) \ge \langle \upsilon, x - x_0 \rangle,
\]
справедливым для всех $x$ из области определения. Отметим, что субградиентые методы широко используются в оптимизации (см., например \cite{gasnikov2017book,nesterov2018lect,polyak1987convex,shor1979nondiff}).

    1. Через центр имеющегося квадрата проводится горизонтальная прямая. На отрезке, высекаемом из квадрата этой прямой, с точностью $\delta$ по аргументу решается задача одномерной оптимизации. Далее при реализации для решения одномерной задачи мы используем метод золотого сечения.

    2. В найденной точке $x_{\delta}$ вычисляется вектор (суб)градиента $\nabla f(x_{\delta})$ и определяется, в сторону какого из двух прямоугольников направлен вектор $\nabla f(x_{\delta})$ и этот прямоугольник исключается из рассмотрения. При этом важно, что достаточно знать именно направление, а не точное значение (суб)градиента.

    3. Через центр оставшегося прямоугольника проводится вертикальная прямая, на отрезке, высекаемом этой прямой в прямоугольнике, также с точностью $\delta$ по аргументу решается задача одномерной оптимизации. В найденной точке вычисляется вектор (суб)градиента функции и определяется, в сторону какого из двух квадратов он направлен. Этот квадрат исключается из рассмотрения.

Если вектор (суб)градиента в точке-приближении решения вспомогательной одномерной задачи нулевой, то процесс можно остановить и выдать указанную точку (она и даст точное решение).

Также возможно, что (суб)градиент в данной точке будет направлен вдоль отрезка, на котором решается вспомогательная одномерная задача и неясно, какую часть оставшейся фигуры отсекать. В таком случае условимся отсекать любую из частей на выбор.

%%%%%%%%%%%%%%%%%%%%%%%%%%%%%%%%%%%%%%%%%%%
%%%%%%%%%%%%%%%%%%%%%%%%%%%%%%%%%%%%%%%%%%%
%%ОБОСНОВАНИЕ ОЦЕНКИ СКОРОСТИ СХОДИМОСТИ%%%
%%%%%%%%%%%%%%%%%%%%%%%%%%%%%%%%%%%%%%%%%%%
%%%%%%%%%%%%%%%%%%%%%%%%%%%%%%%%%%%%%%%%%%%

\section{Обоснование оценки скорости сходимости}

Легко понять, что на каждой итерации метода Ю.Е. Нестерова линейные размеры квадрата уменьшаются вдвое и в какой-то момент оставшийся квадрат будет настолько мал, что в силу условия Липшица значения целевого функционала в его точка будут достаточно близкими друг к другу и можно выбирать любую точку. Мы докажем, что при подходящем выборе точности решения вспомогательных одномерных задач можно гарантировать, что значения целевой функции в точках оставшегося квадрата будут достаточно близки к оптимальному.

    \begin{theorem}

        Пусть на некотором квадрате $\Pi \subset \R^2$ со стороной $R$ задана выпуклая дифференцируемая функция $f: \Pi \rightarrow \R$, для которой при любых $x, y \in \Pi$

        \begin{enumerate}[{1)}]
            \item $\abs{f(x) - f(y)} \leqslant L \|x-y\|$,
            \item $\|\nabla f(x) - \nabla f(y)\| \leqslant M \|x-y\|$.
        \end{enumerate}
        % где $L$ выбирается таким образом, что $L R \sqrt{2}/\varepsilon = 2^x$ для некоторого $x \ge 0$.

        Тогда для всякого $\varepsilon > 0$ после
        \[
            n = \log_2 \frac{2 L R\sqrt{2}}{\varepsilon}
        \]
        шагов метода Ю. Е. Нестерова ($L$ выбрана так, чтобы n было натуральным) при условии, что вспомогательные одномерные задачи решаются с точностью по аргументу
    %     \[
    % \delta = \frac{\varepsilon}{2M R \sqrt{5} \log_2 \frac{2 L R\sqrt{2}}{\varepsilon}},
    %     \]
        \[
        \delta = \frac{\varepsilon}{2 M R \left(\sqrt{2} + \sqrt{5}\right) \left(\displaystyle 1 - \frac{\varepsilon}{L R \sqrt{2}}\right)},
        \]
        получим квадрат $K_n \subset \Pi = K_0$ такой, что
        \[
            f(x) - \min_{x \in \Pi} f(x) \leqslant \varepsilon \quad \forall x \in K_n.
        \]

    \end{theorem}

     \begin{proof}
    \begin{enumerate}[{1.}]
        \item На каждом шаге метода линейные размеры квадрата уменьшаются в 2 раза. Поэтому после $n$ шагов метода длина диагонали оставшегося квадрата $K_n$ будет равна $\displaystyle \frac{R\sqrt{2}}{2^n}$. Из условия Липшица для $f$ имеем
        \[
            f(x) - \min_{x \in K_n} f(x) \leqslant L \|x-\arg\min_{x \in K_n} f(x)\|, \quad \forall x \in K_n.
        \]
        Наибольшее возможное расстояние между точками в полученном квадрате равно длине его диагонали. Это позволяет оценить правую часть неравенства сверху.
        \begin{gather} \label{1}
            f(x) - \min_{x \in K_n} f(x) \leqslant \frac{L R\sqrt{2}}{2^n}, \quad \forall x \in K_n.
        \end{gather}

        Тогда получаем критерий останова метода
        \begin{gather} \label{2}
            \frac{L R \sqrt{2}}{2^n} \leqslant \frac{\varepsilon}{2} \quad\text{или}\quad n \geqslant \log_2 \frac{2 L R\sqrt{2}}{\varepsilon}.
        \end{gather}
        % При этом правая часть последнего неравенства, в силу условия
        % \[
        %     \frac{L R \sqrt{2}}{\varepsilon} = 2^x,\quad x \ge 0
        % \]
        % является натуральным числом, так что критерий останова можно переписать в форме равенства:
        % \[
        %     n = \log_2 \frac{2 L R\sqrt{2}}{\varepsilon}.
        % \]
        % (Очевидно, что возможно подобрать такое $L \geq L_{min}$, что для него выполняется требуемое условие, кроме того, $L$ будет превосходить $L_{min}$ не более чем в два раза. Под обозначением $L_{min}$ понимается константа Липшица функции.)

        \item
        Пусть после некоторого количества итераций метода мы рассматриваем квадрат $P$. На текущей итерации он разделяется на два меньших прямоугольника $Q$ и $Z$, причём на разделяющем отрезке $l$ решается задача одномерной оптимизации $f$.

        Как известно, для её решения $x_0$ можно найти приближение $x_{\delta}$ со любой сколь угодно малой точностью по аргументу $\delta$ (используя такие методы одномерной оптимизации, как, например, метод золотого сечения)
        \[
            \|x_{\delta} - x_0\| \leqslant \delta.
        \]
        
        \item  \begin{enumerate}[{1)}]
            \item Обозначим через $x_0$ решение вспомогательной задачи
            \[
                f(x) \rightarrow \min_{x \in l}.
            \]

        Из липшицевости градиента $f$: $\|\nabla f(x) - \nabla f(y)\| \leqslant M \|x-y\|$ имеем
        \begin{gather} \label{4}
            \|\nabla f(x_{\delta}) - \nabla f(x_0)\| \leqslant M \delta.
        \end{gather}

        Допустим, что $\nabla f (x_{\delta})$ направлен в сторону $Z$, и $x_{\delta}$ не совпадает ни с одним из концов отрезка $l$.  Докажем, что
        \begin{gather} \label{5}
            \min_{x \in Q} f(x) - \min_{x \in P} f(x) \leqslant \frac{M R \delta \sqrt{5}}{2}
        \end{gather}

        Согласно лемме Ферма (не уменьшая общности рассуждений, мы считаем, что решение вспомогательной задачи не совпадает с концом отрезка)
            \[
                x_0 = \arg\min_{x\in l} f(x) \Rightarrow \frac{\partial f}{\partial x_i} = 0.
            \]
            Тогда $\nabla_i f(x_0) = 0$, откуда
            \[
                \left[\begin{gathered}\begin{gathered}
                    \nabla f(x_0) = \overline{0}, \hfill \\
                    \nabla f(x_0)\; \bot\; l. \hfill \\
                \end{gathered}\end{gathered}\right.
            \]

            \item Пусть $\nabla f(x_0) = \overline{0}$. Тогда $x_0 \in l \subset Q$ – точка локального (а вследствие выпуклости $f$ и глобального) минимума функции $f$ на $P$. В таком случае
            \begin{gather} \label{6}
                \min_{x \in Q} f(x) = \min_{x \in P} f(x).
            \end{gather}

            \item Если вектор $\nabla f(x_0) \neq \overline{0}$ и направлен в сторону $Z$, то также верно \eqref{6}, так как производная по всякому направлению, лежащему в $Z$, в таком случае будет неотрицательной.

            \item Если $\nabla f(x_0) \neq \overline{0}$ и направлен в сторону $Q$, то вектора $\nabla f(x_0)$ и $\nabla f (x_{\delta})$ направлены в разные части делимого прямоугольника. Также вектор $\nabla f(x_0)$ перпендикулярен $l$ и поэтому вектора $\nabla f(x_0)$ и $\nabla f (x_{\delta})$ образуют тупой угол.

            В тупоугольном треугольнике, составляемом этими векторами, против большего угла лежит большая сторона, откуда
            \[
                \|\nabla f (x_{\delta})\| \leqslant \|\nabla f(x_{\delta}) - \nabla f(x_0)\| \overset{(4)}{\leqslant} M \delta.
            \]
            Тогда
            \begin{gather} \label{7}
                \|\nabla f(x_{\delta})\| \leqslant M \delta.
            \end{gather}

            Пусть $\displaystyle x_{*} = \arg\min_{x \in Z} f(x)$. Тогда
            \[
                f(x_{*}) - f(x_{\delta}) \geqslant \langle \nabla f(x_{\delta}), x_{*} - x_{\delta} \rangle \Rightarrow
            \] \[
                \Rightarrow f(x_{\delta}) - f(x_{*}) \leqslant \langle -\nabla f(x_{\delta}), x_{*} - x_{\delta} \rangle \leqslant \abs{\langle -\nabla f(x_{\delta}), x_{*} - x_{\delta} \rangle}.
            \]
            По неравенству Коши-Буняковского:
            \[
                f(x_{\delta}) - f(x_{*}) \leqslant \|\nabla f(x_{\delta})\| \cdot |x_{*} - x_{\delta}|.
            \]

            Известно, что $x_{*}$ находится в одном из двух прямоугольников, на которые был поделен исходный квадрат. Диагональ этого прямоугольника (длина которой есть также и наибольшее возможное расстояние между точками в нём) равна $\displaystyle \frac{R\sqrt{5}}{2}$. То есть
            \[
                \|x_{\delta} - x_{*}\| \leqslant \frac{R\sqrt{5}}{2}.
            \]
            Так как в $x_{\delta}$ достигается минимум $f(x)$ на $Q$, а в $x_{*}$ – минимум $f(x)$ на $P$, то верно следующее неравенство
            \[
                \min_{x \in Q} f(x) - \min_{x \in P} f(x) \leqslant \frac{M R \delta \sqrt{5}}{2}.
            \]

            \item Если точка $x_{\delta}$ находится на одном из концов $l$ и градиент в этой точке направлен в сторону $Z$, то все производные по направлениям, отклоняющимся от направления градиента не более, чем на $\displaystyle \frac{\pi}{2}$, будут неотрицательны. Таким образом имеем, что на полуплоскости, отделяемом прямой, проходящей через конец отрезка и содержащим $Z$, значения функции не меньше, чем в самом конце отрезка. Тогда
            \[
                \min_{x \in Q} f(x) = \min_{x \in P} f(x).
            \]
        \end{enumerate}

        \item Пусть на некоторой итерации имеем квадрат $P$ со стороной $R_0$. На первом шаге из рассмотрения исключается один из прямоугольников, на которые был разделён данный квадрат. Длина диагонали оставшегося прямоугольника равна $\displaystyle \frac{R_0 \sqrt{5}}{2}$. При этом ввиду (5) имеем
        \[
        \min_{x \in Q} f(x) - \min_{x \in P} f(x) \le M \delta \cdot \frac{R_0 \sqrt{5}}{2}.
        \]
        На втором шаге пополам делится уже оставшийся прямоугольник, в результате чего получаем квадрат, длина диагонали которого равна $\displaystyle \frac{R_0 \sqrt{2}}{2}$. Соответственно,
        \[
        \min_{x \in W} f(x) - \min_{x \in Q} f(x) \le M \delta \cdot \frac{R_0 \sqrt{2}}{2}.
        \]

                \item На каждой итерации метода линейные размеры квадрата уменьшаются вдвое и поэтому на $i$-ой итерации длина стороны квадрата будет равна $\displaystyle R \cdot \frac{1}{2^{i-1}}$. Каждая же итерация предполагает отсечение сначала прямоугольника, а затем~--- квадрата. Просуммировав оценки накапливающихся погрешностей (см. пункты 3 и 4 выше), получим:
        \[
        \min_{x \in K_n} f(x) - \min_{x \in \Pi} f(x) \le M \delta \cdot \left(1 + \frac{1}{2} + ... + \frac{1}{2^{n-1}}\right) R \cdot \left(\frac{\sqrt{5}}{2} + \frac{\sqrt{2}}{2} \right).
        \]
        С учётом формулы суммы первых $n$ членов геометрической прогрессии
        \[
        \sum_{i=0}^{n-1} \frac{1}{2^i} = \frac{\displaystyle 1-\frac{1}{2^{n-1}}}{\displaystyle 1 - \frac{1}{2}} = 2 \left(1 - \frac{1}{2^{n-1}}\right)
        \]
        получаем
        \begin{equation}\label{75}
        \min_{x \in K_n} f(x) - \min_{x \in \Pi} f(x) \le M R \delta \left(\sqrt{2} + \sqrt{5}\right) \left(1 - \frac{1}{2^{n-1}}\right).
        \end{equation}

        \item
        Если потребовать, чтобы выполнялось неравенство
        \begin{gather} \label{8}
        \min_{x \in K_n} f(x) - \min_{x \in \Pi} f(x) \leqslant \varepsilon / 2,
        \end{gather}
        то ввиду \eqref{2} и \eqref{75} получим
        \[
        M R \delta \left(\sqrt{2} + \sqrt{5}\right) \left(1 - \frac{\varepsilon}{L R \sqrt{2}}\right) \le \frac{\varepsilon}{2},
        \]
        т.е. необходимая точность решения вспомогательных одномерных задач оптимизации:
        \begin{gather} \label{9}
        \delta \leqslant \frac{\varepsilon}{2 M R \left(\sqrt{2} + \sqrt{5}\right) \left(\displaystyle 1 - \frac{\varepsilon}{L R \sqrt{2}}\right)}.
        \end{gather}

        Тогда (3) и (9) означают, что
        \begin{gather} \label{11}
        f(x) - \min_{x \in \Pi} f(x) \leqslant \varepsilon\quad \forall x \in K_n. \qquad
        \end{gather}

        %Тогда точность решения вспомогательных задач по функции

        %\begin{gather} \label{10}
        %\Delta = L \delta = \frac{\varepsilon L}{(4L + 2M R %\sqrt{5}) \log_2 \frac{2 L R\sqrt{2}}{\varepsilon}}.
        %\end{gather}

        % Тогда (2), (3) и (9) означают
        % \begin{gather} \label{11}
        %     f(x) - \min_{x \in \Pi} f(x) \leqslant \varepsilon,\quad \forall x \in K_n. \qquad
        % \end{gather}
    \end{enumerate}
\end{proof}

\begin{remark}
Из доказанного результата по стандартной схеме рассуждений можно доказать и сходимость метода по аргументу в классе сильно выпуклых целевых функционалов (см., например \cite{polyak1987convex}).
\end{remark}

\begin{remark}
    Важно отметить, что схема доказательства основного результата (теорема 1) позволяет учитывать не только погрешность решения вспомогательных одномерных задач, но и погрешность вычисления градиента в точках $x_{\delta}$ (см. пп. 3 и 4 доказательства теоремы 1). Точнее говоря, в точках $x_{\delta}$ вектор градиента функции может вычисляться с точностью $\Delta $: если $\nabla f(x_{\delta})$~--- точная величина градиента, то нам доступен некоторый вектор $v(x_{\delta})$, для которого
    \[
    \|v(x_{\delta}) - \nabla f(x_{\delta}) \| \leqslant \Delta.
    \]

    В том случае, если вектор $v(x_\delta)$ направлен в ту же часть делимого квадрата, что и вектор $\nabla f(x_\delta)$, все рассуждения п. 3 доказательства теоремы 1 остаются неизменными, так как исключается та же часть квадрата, что и при $\Delta = 0$. Если же это не так, то либо вектор $\nabla f(x_\delta)$ направлен в ту часть, где лежит $x^*$ – точное решение задачи, а вектор $v(x_\delta)$ – в иную (очевидно, что в результате будет исключена часть, которой не принадлежит $x^*$, так что рассуждения в п. 3 также остаются справедливыми), либо же вектор $v(x_\delta)$ в ту часть, где лежит $x^*$, а $\nabla f(x_\delta)$ –-- в иную.

    В таком случае имеем треугольник, образованный векторами $\nabla f(x_\delta) - \nabla f(x_0)$ и $v(x_\delta) - \nabla f(x_\delta)$. Из неравенства треугольника следует, что если $\|\nabla f(x_\delta) - \nabla f(x_0)\| \overset{(4)}{\leqslant} M \delta$ и $\|v(x_\delta) - \nabla f(x_\delta)\| \le \Delta$, то
    \[
    \|v(x_\delta) - \nabla f(x_0)\| \le M \delta + \Delta.
    \]

    Векторы $v(x_\delta)$ и $\nabla f(x_0)$ образуют тупоугольный треугольник, причём длина его наибольшей стороны $\|v(x_\delta) - \nabla f(x_0)\| \le M \delta + \Delta$. Следовательно, $\|v(x_\delta)\| \le M \delta + \Delta$. Из треугольника, образуемого векторами $v(x_\delta)$ и $\nabla f(x_\delta)$ ввиду неравенства треугольника, получаем
    \begin{equation}\label{Eqf1}
    \|\nabla f(x_\delta)\| \le \|v(x_\delta)\| + \|v(x_\delta) - \nabla f(x_\delta)\| \le M \delta + 2 \Delta.
    \end{equation}

    Далее, следуя схеме рассуждений доказательства теоремы 1, получим
    \[
    \left(M \delta + 2 \Delta\right) \cdot R \left(\sqrt{2} + \sqrt{5}\right) \left(1 - \frac{\varepsilon}{L R \sqrt{2}}\right) \le \frac{\varepsilon}{2},
    \]
    откуда
    \begin{equation}\label{Eqf2}
    2 \Delta + M\delta \le \frac{\varepsilon}{2R \left(\sqrt{2} + \sqrt{5}\right) \left(\displaystyle 1 - \frac{\varepsilon}{L R \sqrt{2}}\right)}.
    \end{equation}

    Таким образом, метод можно применить при наличии двух типов погрешностей: при решении вспомогательных одномерных задач оптимизации и при вычислении (направлений) градиентов. Эти два типа погрешностей аккумулируют итоговую погрешность. При увеличении одной из погрешности можно уменьшать значение другой так, чтобы сохранялось последнее неравенство, гарантирующее достижение приемлемого качества решения.
    \end{remark}

    \begin{remark}
        Стоит отметить, что в случае, когда вектор-приближение градиента $v(x_{\delta})$ направлен противоположно градиенту в точке решения вспомогательной задачи $\nabla f(x_0)$, в силу полученной выше оценки \eqref{Eqf1}
        $$\|\nabla f(x_\delta)\| \le M \delta + 2 \Delta,$$
        т.е. при достаточно близких к нулю малых величинах $\delta$ и $\Delta$ вектор $\nabla f(x_{\delta})$ можно считать достаточно малым по модулю. C учетом условия \eqref{Eqf2} для суммарной погрешности решения вспомогательных задач
        приходим к следующей оценке значения градиента:
        $$\|\nabla f(x_\delta)\| \le \frac{\varepsilon}{2R \left(\sqrt{2} + \sqrt{5}\right) \left(\displaystyle 1 - \frac{\varepsilon}{L R \sqrt{2}}\right)}.$$
        В таком случае можно считать полученную точку $x_{\delta}$ достаточно хорошим приближением точного решения задачи: используя полученное ранее неравенство
        \[
            f(x_{\delta}) - f(x_{*}) \leqslant \|\nabla f(x_{\delta})\| \cdot |x_{*} - x_{\delta}|,
        \]
        получаем оценку для неточности по функции при выборе такого решения:
        \[
            f(x_{\delta}) - f(x_{*}) \leq \frac{\varepsilon R \sqrt{5}}{4 R \left(\sqrt{2} + \sqrt{5}\right) \left(\displaystyle 1 - \frac{\varepsilon}{L R \sqrt{2}}\right)} \leq \frac{\varepsilon}{4 \left(\displaystyle \sqrt{2} - \frac{\varepsilon}{L R}\right)}.
        \]
        Откуда при дополнения условия теоремы 1 ограничением на $L$:
        \[
            L \geq \frac{\varepsilon}{1.6 R},
        \]
        получаем
        \[
            f(x_{\delta}) - f(x_{*}) \leq \varepsilon.
        \]

       Отправляясь от следующего неравенства для неточного градиента
        \[
            \|\nabla f(x_{\delta})\| \leq \|v(x_{\delta})\| + \Delta,
        \]
        выпишем достаточное условие на вектор $v(x_{\delta})$, при выполнении которого найденная точки будет приемлемой, т.е. гарантированно выполниться \eqref{Eqf1}:
        \[
            \|v(x_{\delta})\| \leq M \delta + \Delta.
        \]
        Таким образом, при использовании данного наблюдения, критерий останова метода можно переписать так:
        \[
            \left(n = \log_2 \frac{2 L R\sqrt{2}}{\varepsilon}\right) \vee \left(\|v(x_{\delta})\| \leq M \delta + \Delta \right).
        \]
    \end{remark}

    Отметим, что при этом можно построить такой пример функции, что для неё даже при небольшой погрешности вычисления градиента может не наблюдаться свойство сходимости метода по аргументу.

    \begin{example}
    Рассмотрим на квадрате $[0; 1]^2$ функцию
    \[
    f(x_1, x_2) = x_1 - 0.001 x_2,
    \]
    минимум которой достигается в точке $(0; 1)$. Градиент данной функции постоянен и равен $\nabla f(x) = (1; -0.001)$, причём вторая координата мала по модулю. Если при некотором $\Delta$ эта компонента будет найдена с погрешностью большей $0.001$, то может оказаться, что она станет положительной и полученный вектор будет направлен в прямоугольник, содержащий точку точного решения задачи, в результате чего он будет исключен из рассмотрения. Из-за этого после удаления прямоугольника, содержащего точное решение задачи на первой итерации, метод не будет сходиться по аргументу.

    Тем не менее, по функции будет достигнуто необходимое качество решения, поскольку ввиду малого коэффициента при $x_2$ слагаемое $0.001 x_2$ не будет значительно влиять на значение функции.
    \end{example}

    \begin{remark}
        Функция, рассмотренная в примере 1, обладает также тем свойством, что константа Липшица её градиента равна 0, что позволяет выбирать в качестве значения $M$ любую сколь угодно малую величину. Поскольку $\delta \sim \frac{1}{M}$ (см. (10)), то в этом случае погрешность решения одномерной подзадачи по аргументу может быть сколь угодно большой. То есть для функций (в том числе и нелинейных), для которых значение $M$ на данном квадрате достаточно мало, может оказаться допустимым вообще не решать вспомогательные задачи на разделяющих отрезках.
    \end{remark}

%%%%%%%%%%%%%%%%%%%%%%%%%%%%%%%%%%%%%%%%
%%%%%%%%%%%%%%%%%%%%%%%%%%%%%%%%%%%%%%%%
%%РАБОТА МЕТОДА НА НЕГЛАДКИХ ФУНКЦИЯХ%%%
%%%%%%%%%%%%%%%%%%%%%%%%%%%%%%%%%%%%%%%%
%%%%%%%%%%%%%%%%%%%%%%%%%%%%%%%%%%%%%%%%

\section{О работе метода Ю.Е. Нестерова для негладких целевых функций}

Итак, мы обосновали сходимость за линейное время метода Ю. Е. Нестерова в классе гладких функций двух переменных с липшицевым градиентом.

\begin{remark}
    На самом деле условие Липшица для градиента в доказательстве теоремы 1 было использовано лишь для точек, лежащих на высекаемых отрезках, вдоль которых рассматривались задачи одномерной минимизации (теорема 1, п. 3 доказательства). Если известна точность $\varepsilon$ для желаемого решения, а также константа Липшица $L$ функционала $f$, то можно определить число итераций $n$ метода Ю. Е. Нестерова. По этому числу $n$ можно построить $2 \cdot (2^n - 1)$ прямых – $2^n - 1$ вертикальных и то же количество горизонтальных, вдоль которых возможно будут решаться вспомогательные задачи.

    Условие Липшица для градиента целевого функционала достаточно потребовать только в точках построенных отрезков. Это, в частности, означает, что вне точек этих отрезков функция даже может иметь точки негладкости.
\end{remark}

\begin{remark}
    Условие липшицевости градиента используется лишь для отрезков, вдоль которых решаются вспомогательные задачи оптимизации, и в силу того, что эти отрезки параллельны осям координат, допустимо использовать условие Липшица не для полного градиента, а лишь для одной из частных производных функции по соответствующей переменной. Если окажется, что по одной из переменных функция линейна, то при решении вспомогательных задач константа Липшица соответствующей частной производной будет равна нулю, что означает, что для некоторых вспомогательных задач $\delta$ может быть сколь угодно велико и эти задачи не требуют решения.
\end{remark}

Естественный интерес представляет вопрос о возможности обобщения теоремы 1 на случай негладкого выпуклого липшицева функционала $f$. Оказывается, что можно привести пример, для которого такая схема будет эффективно работать. Напомним, что для гладких функций $f_1$ и $f_2$ функция

\[
    f(x) = \max\{f_1(x), f_2(x)\}
\]
будет уже негладкой и по теореме Демьянова-Данскина субградиентом $f$ может считаться градиент $\nabla f_1(x)$ при $f_1(x) \ge f_2(x)$ и наоборот $\nabla f_2(x)$ при $f_2(x) \ge f_1(x)$.

Приведем результаты работы методов для функции $$f(x_1, x_2) = \max\{2x_1 + 6; -x_1; x_2 + 3; -2x_2\}.$$

Критерий остановки для каждого из методов --- достижение точности $\varepsilon = 5 \times 10^{-3}$ по функции в сравнении с известным точным решением). Для метода Ю.Е. Нестерова время работы составляет около 2 ms, для градиентного спуска~--- около 5 ms, для метода эллипсоидов~--- значительно дольше.

Тем не менее, возможен простой пример негладкой выпуклой функции, для которой рассматриваемый метод не будет сходиться.

\begin{example}
    Рассмотрим на квадрате с вершинами в точках $(0; 0), (0; 1), (1; 0)$ и $(1; 1)$ функцию
    \[
        f(x_1, x_2) = \abs{x_1 - x_2} + 0.9x_1.
    \]
    В качестве первой разделяющей прямой выберем $x_2 = 1/2$. Тогда при $x_1 < 1/2$ имеем $f\left(x_1, x_2\right) = x_2 - 0.1x_1$, а при $x_1 \ge 1/2$ получаем $f\left(x_1, x_2\right) = 1.1x_1 - x_2$.
    Поэтому
    \[
        \nabla f\left(x_1, 1/2\right) = \begin{cases} (-0.1; 1), & \mbox{если } x_1 < 1/2, \\ (1.1; -1), & \mbox{если } x_1 > 1/2. \end{cases}
    \]
    Иными словами, в зависимости от выбора точки приближения решения (и даже при точном решении) вспомогательной задачи (суб)градиент может быть направлен в любую из полуплоскостей относительно $x_2 = 1/2$ (отметим, что $\displaystyle \min_{0 \le x_1 \le 1} f\left(x_1, 1/2\right) = f\left(1/2, 1/2\right)$). Вместе с тем при отбрасывании части квадрата, для которой $x_2 < 1/2$, наименьшее значение целевой функции $f$ в оставшейся части будет $f\left(1/2, 1/2\right) = 9/20 > f(0, 0) = 0$.
\end{example}

%%%%%%%%%%%%%%%%%%%
%%%%%%%%%%%%%%%%%%%
%%ДВА ОГРАНИЧЕНИЯ%%
%%%%%%%%%%%%%%%%%%%
%%%%%%%%%%%%%%%%%%%

\section{О применении метода Ю.Е. Нестерова для решения задач выпуклого программирования с двумя ограничениями}

\subsection{Постановка задачи}

Рассмотрим задачи следующего вида

\begin{equation}\label{st_eq001}
f(x)\rightarrow\min,\;x\in Q,\;g_1(x)\leqslant0,\;g_2(x)\leqslant0,
\end{equation}
где $f,g_1$ и $g_2$~--- выпуклые функционалы, $Q$~--- выпуклый компакт в $\mathbb{R}^n$. Допустим также, что $f$~--- $\mu$-сильно выпукла и $\forall i=\overline{1,2}$
\begin{equation}\label{st_eq002}
|g_i(x)-g_i(y)|\leqslant M_i\|x-y\|
\end{equation}
при всяких $x,y\in Q$ для фиксированных $\mu_{1,2}>0$.

Тогда двойственная задача к \ref{st_eq001} будет имеет такой вид
\begin{equation}\label{st_eq003}
\varphi(\lambda_1,\lambda_2)=\min_{x\in Q}\fbr{f(x)+\lambda_1g_1(x)+\lambda_2g_2(x)}\rightarrow\max_{\lambda_1+\lambda_2\leqslant\Omega_{\lambda}},
\end{equation}
где $\Omega_{\lambda}$ определяется условием Слейтера.

Здесь функция $\varphi$ зависит от двух двойственных переменных $\lambda^1$ и $\lambda^2$, удовлетворяет условию Липшица и имеет липшицев градиент:
\[
    \abs{\varphi(\lambda_1) - \varphi(\lambda_2)} \leq C \|\lambda_1 - \lambda_2\|,
\]
\[
    \|\nabla \varphi(\lambda_1) - \nabla \varphi(\lambda)\| \leq \frac{M^2}{\mu} \|\lambda_1 - \lambda_2\|,
\]
где $C = \max_{x \in Q} \left\{g_1(x), g_2(x)\right\}$, $M = \sqrt{M_1^2 + M_2^2}$.

Если для прямой задачи выполняется условие Слейтера, можно получить ограничение следующего вида для точки-решения двойственной задачи:
\[
    \lambda^*_1 + \lambda^*_2 \leq A.
\]

Таким образом, решение задачи сводится к оптимизации функции $\varphi$ двух переменных на прямоугольном треугольнике $$\Omega_{\lambda}=\Lambda_A = \{\lambda \in \mathbb{R}^2_{+} | \lambda_1 + \lambda_2 \leq A\},$$ 
катеты которого равны $A$ и лежат на осях координат. Для решения данной задачи можно использовать и модифицированный метод Ю. Е. Нестерова для прямоугольного треугольника, который описан далее. Можно рассматривать 2-мерную двойственную задачу и на квадрате, содержащем указанный прямоугольный треугольник.

\subsection{Приложения к задачам выпуклого программирования с двумя ограничениями: адаптивный критерий остановки метода}

Если для некоторого $x_{\delta}(\lambda)$ и $\lambda=(\lambda_1,\lambda_2)$ верно
\begin{equation}\label{st_eq004}
|\lambda_1g_1(x_{\delta}(\lambda))+\lambda_2g_2(x_{\delta}(\lambda))|\leqslant\varepsilon,
\end{equation}
где $x_{\delta}(\lambda)$~--- приближённое решение задачи
\begin{equation}\label{st_eq005}
x(\lambda)=arg\min_{x\in Q}\fbr{f(x)+\lambda_1g_1(x)+\lambda_2g_2(x)}
\end{equation}
с точностью $\delta$ по функции, то
$$f(x_{\delta}(\lambda))-f(x_*)\leqslant\varepsilon+\delta.$$

Действительно, пусть $\lambda_*=(\lambda_1^*,\lambda_2^*)$~--- решение задачи \eqref{st_eq003} и\\ $x_*=x(\lambda_*)$. Тогда
$$f(x_{\delta}(\lambda))+\lambda_1g_1(x_{\delta}(\lambda))+\lambda_2g_2(x_{\delta}(\lambda))\leqslant$$
$$f(x(\lambda))+\lambda_1g(x(\lambda))+\lambda_2g(x(\lambda))+\delta=\varphi(\lambda_1,\lambda_2)+\delta\leqslant$$
$$\varphi(\lambda_*)+\delta=f(x_*)+\lambda_1^*g(x_*)+\lambda_2^*g(x_*)+\delta\leqslant f(x_*)+\delta,$$
откуда
$$f(x_{\delta}(\lambda))-f(x_*)\leqslant-(\lambda_1g_1(x_{\delta}(\lambda))+\lambda_2g_2(x_{\delta}(\lambda)))+\delta\leqslant\varepsilon+\delta,$$
что и требовалось.

Поэтому задачу \eqref{st_eq001} можно сводить к \eqref{st_eq003} и решать \eqref{st_eq003} методом Ю.\,Е.\,Нестерова на подходящем квадрате (или прямоугольном треугольнике) с критерием остановки \eqref{st_eq004} при $\delta \leqslant \varepsilon$. Если применить подходы, описанные в замечаниях 2 и 3, то всякий раз после удаления прямоугольников в текущем промежутке локализации двойственной переменной будет сохраняться оптимальное значение двойственных множителей. Этот промежуток локализации будет уменьшаться. В какой-то момент возмущенный градиент двойственной функции в точке-приближении вспомогательной одномерной задачи будет достаточно малым для выполнения критерия остановки \eqref{st_eq004}. Если решать задачи \eqref{st_eq005} с точностью $\varepsilon$ по функции и аргументу, то можно по схеме рассуждений \cite{Stonyakin_MOTOR_2019} доказать для гладких функционалов $f$, $g_1$ и $g_2$ оценку скорости сходимости $O\br{\log^3\frac{1}{\varepsilon}}$.
%%%%%%%%%%%%%%%
%%%%%%%%%%%%%%%
%%ТРЕУГОЛЬНИК%%
%%%%%%%%%%%%%%%
%%%%%%%%%%%%%%%

\section{Модификация метода для оптимизации на прямоугольном треугольнике}

Приведём одну из возможных модификаций метода Ю.Е. Нестерова для оптимизации функции на прямоугольном треугольнике.

    1. На первом шаге итерации метода проводится разделющий отрезок, соединяющий середины одного из катетов и гипотенузы. С помощью одного из методов одномерной оптимизации (например, метода золотого сечения) решается вспомогательная задача минимизации на данном отрезке.

    2. Далее, в точке-решении вычисляется градиент $\nabla \varphi(x)$, после чего отсекается та часть треугольника, в которую он направлен. В случае, если после отсечения остаётся в два раза меньший гомотетичный исходному треугольник, метод переходит к следущей итерации.

    3. Иначе, на втором шаге итерации, аналогичным образом отсекается одна из частей оставшейся трапеции, на которые её делит отрезок, соединяющий середины гипотенузы и другого катета исходного треугольника. Если в результате этого остаётся треугольник, то метод переходит к следующей итерации.

    4. Если же остаётся квадрат, то для дальнейшей оптимизации используется описанный ранее метод Ю.Е. Нестерова для оптимизации функции на квадрате.

Имеет место следующий результат (общая схема доказательства аналогична доказательству теоремы 1).

\begin{theorem}
    Пусть на прямоугольном треугольнике $\Lambda_A \in \mathbb{R}^2$ со стороной $A$ задана выпуклая дифференцируемая функция $f: \Lambda_A \rightarrow \mathbb{R}$, для которой при любых $x, y \in \Lambda_A$
    \begin{itemize}
        \item[1)] $|f(x) - f(y)| \leq L \|x - y\|$,
        \item[2)] $\|\nabla f(x) - \nabla f(y)\| \leq M \|x - y\|$.
    \end{itemize}
    Тогда для всякого $\varepsilon > 0$ после $i$ итераций предложенной модификации метода и
    $$
        n - i = \log_2 \frac{2 L A \sqrt{2}}{\varepsilon} - i
    $$
    итераций метода Ю.Е. Нестерова ($L$ подбирается так, чтобы $n-i$ было целым) для оптимизации на квадрате при условии, что вспомогательные задачи решаются с точностью по аргументу
    $$
        \delta \leq \frac{\varepsilon}{2 M A (\sqrt{2} + \sqrt{5}) \left(\displaystyle1 - \frac{\varepsilon}{L A \sqrt{2}}\right)}
    $$
    получим квадрат $K_n \subset \Lambda_A$ такой, что
    $$
        f(x) - \min_{x \in \Lambda_a} f(x) \leq \varepsilon \quad x \in K_n.
    $$
\end{theorem}

%%%%%%%%%%%%%%%%%
%%%%%%%%%%%%%%%%%
%%ЭКСПЕРИМЕНТЫ%%%
%%%%%%%%%%%%%%%%%
%%%%%%%%%%%%%%%%%

\section{Эксперименты}
    Основное преимущество метода Ю.Е. Нестерова состоит в том, что данный метод является неполноградиентным, то есть не требует для своей работы вычисления полного значения градиента функции на каждой итерации. При этом некоторым неодостатком может быть необходимость знать константы $L$ и $M$ (или их верхние оценки) для целевой функции.  В то же время работа, например, метода эллипсоидов не зависит от уровня гладкости целевой функции. Также может оказаться не вполне удобным, что необходимая точность решения вспомогательных одномерных задач зависит от $\varepsilon$, что требует полного перезапуска метода при изменении точности решения задачи.

    Приведем результаты работы методов для функции $f(x_1, x_2) = (x_1-1)^2 + x_2^4$ (критерий останова для каждого из методов~–-- достижение необходимой точности по функции в сравнении с известным точным решением). Для метода Ю.Е. Нестерова время работы составляет около 3 ms, для градиентного спуска~--- более 7 ms, для метода эллипсоидов~--- около 4 ms.

    Рассмотрим гладкую выпуклую функцию $f(x_1, x_2) = (x_1-1)^2 + x_2^4$ на квадрате с вершинами в точках $(-3; -3), (-3; 1), (1; 1), (1; -3)$. В данном эксперименте критерием остановки методов будет достижение необходимой точности по функции $\varepsilon = 5 \times 10^{-3}$ в сравнении c минимумом функции на квадрате, который достигается в точке $(1; 0)$.

    Эксперименты показали, что в сравнении с другими рассматриваемыми методами (градиентным спуском и методом эллипсоидов), метод Ю.Е. Нестерова работает быстрее, несмотря на большее количество выполненных итераций. Это достигается за счёт того, что другие рассмотренные полноградиентные методы тратят большее время на вычисление полного значения градиента функции. При этом рассматриваемый в работе метод Ю.Е. Нестерова предполагает для своей работы лишь вычисление направления (суб)градиента, что, вообще говоря, может удешевлять стоимость итерации.

    Приведем результаты работы методов для функции $f(x_1,x_2) = (x_1 + 1)^2 + x_2^2 - x_1 + e^{x_1} + e^{x_2+1}$ (остановка для каждого из методов определяется в соответствии с теоретическими оценками (см., например,  \cite{polyak1987convex})). Для метода Ю.Е. Нестерова время работы составляет около 8 ms, для метода эллипсоидов~--- более 15 ms.

%%%%%%%%%%%%%%%%%%%%%%%%%%%%%%%%%%%%%%
%%%%%%%%%%%%%%%%%%%%%%%%%%%%%%%%%%%%%%
%%ОБСУЖДЕНИЕ ПОЛУЧЕННЫХ РЕЗУЛЬТАТОВ%%%
%%%%%%%%%%%%%%%%%%%%%%%%%%%%%%%%%%%%%%
%%%%%%%%%%%%%%%%%%%%%%%%%%%%%%%%%%%%%%

\section{Заключение}

В работе получена оценка скорости сходимости одного неполноградиентного метода минимизации выпуклой липшицевой функции двух переменных на квадрате.Как показали результаты экспериментов, такой метод может работать существенно быстрее метода эллипсоидов.

Основной результат (теорема 1) преполагает помимо условия Липшица для целевой функции ещё и условие липшицевости градиента. Однако как показано далее, это требование можно существенно ослабить, потребовав его выполнение на отдельных отрезках (замечание 5). Очевидно, что если, например, точек негладкости целевого функционала конечно, то незначительным сдвигом разделяющих (на итерациях метода) квадрат отрезков можно исключить эти точки негладкости из рассмотрения. Однако нами предложен также и простой пример негладкой функции, для которой метод не сходится (пример 2).

Весьма интересной была бы задача исследовать эффективность предложенного подхода к условным задачам многомерной оптимизации с небольшим числом ограничений (которые за счёт операции максимизации можно свести к двум). Заметим, что выполнение условия Слейтера позволяет компактифицировать значения двойственных переменных (см., например раздел 4 из \cite{gasnikov2017book}) и по сути рассматривать в качестве двойственной задачи двумерную задачу оптимизации (на квадрате или, несколько модифицировав предложенную схему, на прямоугольном треугольнике). В пункте 5 описан общий подход к таким задачам на базе метода Ю.Е. Нестерова для 2-мерной двойственной задачи.

Все численные эксперименты были выполнены с помощью с помощью программного обеспечения Anaconda 5.3.1 Python 3.6 version \cite{anaconda} на компьютере с 2-ядерным процессором Intel Core i5-5250U с тактовой частотой 1,6 ГГц на каждое ядро. ОЗУ компьютера составляла 8 Гб.

\end{document}